\documentclass[11pt]{article}
\newcommand{\comm}[1]{}

\usepackage[T1]{fontenc}
\usepackage[utf8]{inputenc}
\usepackage{color}

\usepackage{amssymb}
\usepackage{amsthm}
\usepackage[centertags]{amsmath}
\usepackage{amsfonts}
\usepackage{amssymb}
\usepackage{amsthm}
\usepackage{epsfig}
\usepackage{setspace}
\usepackage{ae} 
\usepackage{eucal} 
\def\citet{\cite}

   \csname @addtoreset\endcsname{equation}{section}

\newcounter{thanksnum}
\def\thanksnumber#1
{\setcounter{thanksnum}{\value{footnote}}\setcounter{footnote}{#1}%
                     \addtocounter{footnote}{-1}\footnotemark
                     \setcounter{footnote}{\value{thanksnum}}}
\newtheorem{theorem}{Theorem}[section]
\newtheorem{lemma}{Lemma}[section]

\newtheorem{proposition}{Proposition}[section]
\newtheorem{corollary}{Corollary}[section]

\newtheorem{definition}{Definition}[section]
\newtheorem{remark}{Remark}[section]

\newtheorem{example}{Example}[section]

\def\e{\varepsilon}

\def\defi{\stackrel{{\scriptscriptstyle \Delta}}{=}}

\def\a{\alpha}
\def\d{\delta}

\def\O{\Omega}

\def\Y{{\cal Y}}
\def\F{{\cal F}}
\def\w{\widehat}
\def\Ind{{\mathbb{I}}}

\def\const{{\rm const\,}}

\def\Var{{\rm Var\,}}

\def\R{{\bf R}}
\def\E{{\bf E}}
\def\P{{\bf P}}

\def\Z{{\cal Z}}

\def\H{{\cal H}}
\def\T{{\cal T}}

\def\L{L}

\def\s{\delta}
\def\g{\gamma}

\def\X{{\cal X}}
\def\t{\theta}
\def\oo{\bar}
\def\s{\sigma}

\def\G{\Gamma}
\def\GG{{\cal G}}

 \def\V{{\cal V}}

\def\L{{\cal L}}

\newcommand{\be}{\begin{equation}}
\newcommand{\ee}{\end{equation}}
\newcommand{\bd}{\begin{displaymath}}
\newcommand{\ed}{\end{displaymath}}
\newcommand{\ba}{\begin{array}{ll}}
\newcommand{\ea}{\end{array}}
\newcommand{\baa}{\begin{eqnarray}}
\newcommand{\eaa}{\end{eqnarray}}
\newcommand{\baaa}{\begin{eqnarray*}}
\newcommand{\eaaa}{\end{eqnarray*}}   \font\sm=cmr10

\def\H{{\cal H}}

\def\Ts{s}
\def\DD{{\cal D}}

\def\TT{\Theta}
\def\dd{d_{\!{\scriptscriptstyle \bf-^{\vphantom{x}}}\!\!}}
\def\dd{{d_{\scriptscriptstyle F}}}
\def\PC{{{\scriptscriptstyle PC}}}
\def\RS{{{\scriptscriptstyle RS}}}

\def\cH{}
\def\ddH{d_H}

\date{Submitted: October 5, 2015. Revised: April 19, 2020} 
\title{On the  fractional stochastic  integration for random non-smooth integrands}
\author{
Nikolai Dokuchaev\\
 {\sm  School of Electrical Engineering, Computing and Mathematical Sciences, Curtin
University} }
\begin{document}
\maketitle
\begin{abstract}
The paper suggests a way of stochastic integration of random integrands with respect to
fractional Brownian motion with the Hurst parameter $H> 1/2$.   The integral is defined  initially on the processes  that are "piecewise" predictable on a short horizon.  Then  the integral is extended on a wide class of square integrable adapted  random processes. This class is described via a mild
restriction on the growth rate  of the conditional mean square error  for the forecast on an arbitrarily short horizon given current observations; differentiability of  H\"older property of any kind  or degree is not required for the integrand.
The suggested integration can be interpreted as foresighted  integration for integrands
featuring corresponding restrictions on the forecasting error.
This integration is based on It\^o's integration and does not involve
Malliavin calculus or Wick products.
In addition, it is shown that  these  stochastic integrals depend continuously on  $H$ at  $H=1/2+0$.

\end{abstract}

 \section{Introduction}
The paper considers stochastic integration of random integrands with respect to fractional Brownian motion. These integrals
can be defined  using different approaches; see
review and discussion in  \cite{Alos,AN,Carm,Ci,DecU,Dec0,Dec,Duncan,Fe,Ma,PT1,PT2,Pri,Zahle}.
This integration has many applications in statistical modelling, especially for quantitative finance; see e.g. \cite{BKiev,B11,BPS,Bi,Bj,C,G,H,H2,HZhou,Ma,R,Shi,Sal}. Special statistical inference methods developed for these models; see e.g. \cite{CNS,Es,GN,Mu}.

Naturally, the integral can be defined as a Riemann sum for piecewise constant in time integrands; the problem is an extension on more general classes of integrands.
There is a special approach base on the so-called the Wick product rather than Riemann sums; see,  e.g. \cite{BKiev,B11,BPS,Bj,Duncan}. This approach allows integrands of quite general type but the
features the Wick product makes the corresponding integrals  quite distinctive from the integrals based on the Riemann sums.

Currently, stochastic integrals  with respect to the fractional Brownian motion $B_H$ with a Hurst parameter $H\in (1/2,1)$
 are defined for  random integrands in the following cases.
 \begin{enumerate}
\item The integral is defined for the integrands that are pathwise H\"older  with index $p> 1-H$; see, e.g., Theorem 21 in \cite{Fe} and  \cite{Ci,Zahle}.
\item  The integral is defined pathwise for  integrands that has $q$-bounded variation with  $q<1/(1-H)$; see, e.g., \cite{Young,Ber}.
\item The integral is defined as  a Skorohod integral for integrands $\g$ \index{$\g\in {\rm Dom}\d$} such that
$\nabla \g$ is $L_p$-integrable for $p>(1/2-H)^{-1}$, where $\nabla$ is the Gross-Sobolev derivative   ( Theorem 3.6 \cite{Dec0} (2003) or Theorem 6.2 \cite{Dec}). This approach is based on anticipating integrals
(see, e.g., \cite{BKiev,Bj,DecU,Duncan}, and review in \cite{Dec}). It can be noted that this requires certain differentiability of the integrand in the sense of  existence of
$\nabla g$ or the fractional derivative  \cite{Alos}.
 \end{enumerate}
 \index{Comment: \cite{Dec}, p.11, remark 3.3 - integral in \cite{DecU} is actually a Wick product in Duncan }
We exclude from this list the integrals based on the Wick product and integrals for piecewise constant integrands.


In this paper, we readdress stochastic integration of random integrands  with respect to fractional Brownian motion. We suggests an integration scheme allowing  to extend the class of admissible random integrands known in the literature.
In particular, we show that stochastic integral with respect to the fractional Brownian motion $B_H$ with $H\in (1/2,1)$
is well defined  on a wide class of $L_2$-integrable processes with a mild  restriction on the growth rate
for conditional variance for a  short term forecast. It is not required that the integrands $g$  satisfy H\"older
condition, or have finite  $p$-variation, or $\nabla \g$  is $L_p$- integrable, or a fractional derivative exists.  The description of this class does not require to use
Malliavin calculus as in \cite{Dec0, Dec} and does not use any kind of derivatives.

We use a modification of the classical Riemann sums.
 Instead of the standard  extension of the Riemann sums
from the set of piecewise constant integrands,
we used an extension of different sums from
  processes being  "piecewise predictable" on a short horizon  that are not necessarily piecewise constant. More precisely, these
 integrands are  adapted to the filtration
generated by the observations being frozen at grid time points.
In other words, this  "piecewise predictable"  class includes all integrands that  are predictable without error on a fixed  time horizon that can be arbitrarily short.
The corresponding stochastic integral is represented via sums of integrals of two different types: one type is a standard  It\^o's integral, and another
 type is a  Lebesgue integral for random integrands.

In the second step, we extended this integral on a wide class of $L_2$-integrable processes (Theorem \ref{ThM} below); the resulting  integrals is denoted as $\int \cdot \dd B_H$
The corresponding condition allows a simple  formulation that does not require Malliavin calculus  used in \cite{Dec0,Dec}.
This theorem implies prior estimates  of the stochastic integral via a norm of  a random integrand (Corollary \ref{corrI}).

Furthermore, it is shown that the stochastic integrals depend continuously on  $H$ at  $H=1/2+0$ under some additional mild restrictions on the growth rate
for the conditional variance
of the future values given current observations (Theorem \ref{ThL} below).

The paper is organized as  follows. Section \ref{SecD} presents some definitions.
In Section \ref{SecM}, we present the definition  of the new type of integral and
some convergence results and prior estimates. In Section \ref{SecC}, we show some continuity
of the new integral with respect to a variable Hurst parameter.
The proofs are given in Section \ref{SecP}.
\section{Some definitions}\label{SecD}
 We are
given a probability space $(\Omega,\F,\P)$, where $\Omega$ is
a set of elementary events, $\F$ is a complete $\s$-algebra of
events, and $\P$ is a probability measure.

 We  assume that $\{B_H(t)\}_{t\in\R}$ is a fractional Brownian motion
 with the Hurst parameter $H\in  (1/2,1)$ defined as described in \cite{Ma}
 such that $B_H(0)=0$ and
 \baa
B_H(t) =\int_{-\infty}^t f(t,r)dB(r),
 \label{BBH}
 \eaa
 where $t\ge 0$ and \baa
f(t,r)\defi c_H(t-r)^{H-1/2}\Ind_{r\ge 0}+ c_H((t-r)^{H-1/2}-(-r)^{H-1/2})\Ind_{r<0}.
\eaa
 Here
 $c_H=1/\G (H+1/2)$, $\G$ is the Gamma function, $\Ind$ is the indicator function,  and
$\{B(t)\}_{t\in\R}$ is  a standard Brownian motion such that $B(0)=0$; we denote by $\int\cdot dB$ the standard It\^o's integration.

Let $d_H\defi c_H(H-1/2)$.

For ${\rm T}>0$,   $\tau\in [0,{\rm T}]$ and  $g\in L_2(0,{\rm T})$, set \baa
G_H(\tau,{\rm T},g)\defi  d_H \int_{\tau}^{\rm T}(t-\tau)^{H-3/2}g(t)dt.
\label{GH}\eaa
By the property of the Riemann--Liouville integral,   there exists $c>0$ such that
\baa
\|G_H(\cdot,{\rm T},g)\|_{L_2(s,{\rm T})}\le c\|g\|_{L_2(s,{\rm T})}.
\label{estG}
\eaa  It can be noted that this $c$
is independent on $H\in(1/2,1)$.

\par
Let $\{\GG_t\}$ be the filtration generated by the process $B(t)$.

Let $T>0$  be given.

Let  $\L_{22}$ be the linear normed space formed as the completion in $L_2$-norm of the set of all
$\GG_t$-adapted bounded measurable processes $\g(t)$, $t\in[0,T]$, with the norm
$\|\g\|_{\L_{22}}=\left(\E\int_0^T\g(t)^2dt)\right)^{1/2}$.

For $\e>0$, let $\X_\e$  be the set of all  $\g\in\L_{22}$
such that there exists an integer $n>0$ and a set of non-random times
$\T=\{T_k\}_{k=1}^n\subset [0,T]$, where  $T_0=0$, $T_n=T$,
and $0<T_{k+1}-T_k\le \e$, such that $\g(t)$ is $\GG_{T_k}$-measurable for $t\in[T_k,T_{k+1})$.

In particular, the set $\X_\e$ includes all $\g\in\L_{22}$ such that $\g(t)$ is $\GG_{t-\e}$-measurable for all $t\in[0,T]$.

Let $\X_{}\defi \cup_{\e>0}\X_\e$.

For the brevity, we sometimes denote $L_p(\O,\GG_T,\P)$ by $L_p(\O)$, $p\ge 1$.
\index{Let   $B_{H,n}\in\X_{}_{\e,2,\PC}$ be defined as $B_{H,n}(t)=B_H(T_k)$, $t\in[T_k,T_{k+1})$.}

Let $\X_{\e,\PC}$  be the set of all  $\g\in\L_{22}$
such that there exists an integer $n>0$ and a set of non-random times
$\T=\{T_k\}_{k=1}^n\subset [0,T]$, where $n>0$ is an integer, $T_0=0$, $T_n=T$,
and $T_{k+1}-T_k\ge \e$, such that $\g(t)=\g(T_k)$ for $t\in[T_k,T_{k+1})$.
\section{The  main result: integration  for random integrands}\label{SecM}
For any $\g\in\X_{\e,\PC}$, it is naturally to define
the stochastic integral  with respect to $B_H$ in $L_1(\Omega,\GG_{T},\P)$
 as the Riemann sum
 \baaa
 \sum_{k=0}^n\g(T_k)(B_H(T_{k+1})-B_H(T_{k})).
 \eaaa
 If  $\g\in\L_{22}$ is such that this  sum has
  a limit in probability  as $n\to +\infty$,
  and this limit is independent on the choice of  $\{T_k^n\}_{k=1}^n$, then we call this limit
  the integral $\int_0^T\g(t)d_{\RS}B_H(t)$.

The classes of  admissible deterministic integrands $\g$ are known; see, e.g. \cite{PT1,PT2}.  However, there are some difficulties with identifying classes of admissible random $\g$.  The present paper suggests a modification of the stochastic integral based
 on the extension from $\X_{}$, i.e. from the set of random functions  that are not necessarily  piecewise constant but rather "piecewise predictable". This modification will allow to
establish  a new extended class of random integrands that are not necessarily  "piecewise predictable".

\subsubsection*{The case of of non-random integrands}

As the first step, let us construct  a stochastic  integral over the time interval $[s,T]$ for $\GG_s$-measurable integrands  $\g\in L_2(\O,\GG_s,\P,L_2(s,T))$.
 These integrands can be regarded as non-random on the conditional probability space given $\GG_s$.

By  (\ref{BBH}), we have that
\baaa
 B_H(t)=W_H(t)+R_H(t),
 \eaaa
where $t>s$,
 \baaa
 W_H(t)=\cH \int_{s}^t f(t,r)dB(r),\qquad
R_H(t)=\cH   \int_{-\infty}^{s}f(t,r)dB(r).
\label{BBHa}\eaaa
The processes $W_H(t)$ and $R_H(t)$ are independent Gaussian processes with zero mean. In addition, the process $W_H$ is $\{\GG_t\}$-adapted,  $R_H(t)$ is $\GG_s$-measurable for all $t>s$, and   $W_H(t)$ is  independent on $\GG_s$  for all $t>s$.

To define integration with respect to $dB_H$ for $\GG_s$-measurable integrands
 $\g\in L_2(\O,\GG_s,\P,L_2(s,T))$ we define
integration with respect to $W_H$ and $R_H$ separately.

First, it can be noted that if we had $f_t'(t,\cdot)\in L_2(s,t)$ then
integration with respect to $W_H$ would be straightforward, since
 we would be able to find the It\^o's differential $dW_H(t)$ as
\baa
&&f(t,t) dB(t)+\int_0^t f'_t(t,r)dB(r)\cdot dt= 0\cdot dB(t)+\int_0^t f'_t(t,r)dB(r)\cdot dt,
\label{dW}
\eaa
which would allow us to accept  $\int_{s}^T \g(t) \left[\int_0^t f'_t(t,r)dB(r)\right]dt$  as
 $\int_{s}^T\g(t)d W_H(t)$.
However, the expression (\ref{dW}) cannot be regarded as an It\^o's differential, since
$f_t'(t,\cdot)\notin L_2(s,t)$. Nevertheless, we will be using
a modification of this version of the integral with respect to $W_H$  amended
with some approximations to overcome insufficient integrability of $f'_t(t,\cdot)$.

For $\e>0$, let
 \baaa
 W_{H,\e}(t)=\cH \int_{s}^t f(t,r-\e)dB(r).
\label{We}\eaaa
In this case, there  exists a usual  It\^o's differential
\baaa
dW_{H,\e}(t)=f(t,t-\e) dB(t)+\int_0^t f'_t(t,r-\e)dB(r)\cdot dt.
\label{dWe}
\eaaa
representing a "regularized" approximation of the right hand part of (\ref{dW}).  
\begin{proposition}\label{propW}
For any $\g\in L_2(\O,\GG_s,\P,L_2(s,T))$,
 \baaa
\lim_{\e\to 0}\int_{s}^T\g(t)dW_{H,\e}(t) =\int_{s}^TG_H(\tau,T,\g)d B(\tau);
\eaaa
the limit holds in  $L_2(\O,\GG_T,\P)$.
\end{proposition}
This result justifies the following definition.
\begin{definition} \label{defdw}
 We regard the limit in Definition \ref{defdw} as the stochastic integral with respect to $W_H$,
 and we denote it as $\int_{s}^T\g(t)\dd W_H(t)$, i.e.
 \baaa
\int_{s}^T\g(t)\dd W_H(t)\defi \int_{s}^TG_H(\tau,T,\g)d B(\tau).
\eaaa
\end{definition}
It appears that this choice for the case of non-random integrands leads to a new version of a stochastic integral for random integrands constructed below.
\begin{proposition}\label{propM}
 \begin{enumerate}
\item
 $R_H(t)$ is $\GG_s$-measurable for all $t>s$ and differentiable in $t>s$
 in the sense  that
 \baa
\lim_{\d\to 0} \E\left|\frac{R_H(t+\d)-R_H(t)}{\d} -\DD R_H(t)\right|=0, \label{deflim}
\eaa
where
\baaa
\DD R_H(t)\defi \cH   \int_{-\infty}^{\Ts}f'_t(t,q)dB(q).
\eaaa
The process   $\DD R_H$ is such that
 \subitem(a)  $\DD R_H(t)$ is $\GG_s$-measurable for all $t>s$;
\subitem(b) for any $t>s$,
\baa &&\E \DD R_H(t)^2=\frac{d_H^2}{2-2H}(t-\Ts)^{2H-2},\\
&&\E\int_{s}^t \DD R_H(r)^2dr=\frac{c_Hd_H}{2(2-2H)}(t-\Ts)^{2H-1}.\label{EDR}
\eaa
 \end{enumerate}
\end{proposition}

\begin{definition} \label{defWR}
For $s\in [0,T)$ and $\g\in L_2(\O,\GG_s,\P,L_2(s,T))$, we define the integral
\baaa
&&\int_{s}^T\g(t)\dd B_H(t)\defi  \int_{s}^T\g(t)\dd W_H(t)+\int_{s}^T\g(t)\DD R_H(t)dt\\
&&=\int_{s}^TG_H(\tau,T,\g)d B(\tau)+\int_{s}^T\g(t)\DD R_H(t)dt.
\eaaa
The first integral in the sum above is described in Definition \ref{defdw}, and the second one
is a pathwise Lebesgue integral on $[s,T]$. The sum belongs  to $L_1(\O,\GG_T,\P)$ thanks to Propositions    \ref{propW} and \ref{propM}.
\end{definition}

\begin{proposition}\label{propWR} Under the assumptions and notations of
Definition \ref{defWR},
\baaa &&\E\left|\int_{s}^T\g(t)\dd W_H(t)\right|^2\le c\E\int_s^T \g(t)^2dt,
\label{estW} \\
&&\E\left|\int_{s}^T \g(t)\DD R_H(t)dt\right|\le c \left(\E\int_{s}^{T}\g(t)^2dt\right)^{1/2},\label{estR}\\
&&\E\left|\int_{s}^T\g(t)\dd B_H(t)\right|
\le c\left(\E\int_{s}^{T}\g(t)^2dt\right)^{1/2},
\label{estRW}\eaaa
for some $c=c(H,T)>0$.
\end{proposition}
\begin{remark} For the purposes of the proofs below, we need stronger estimates 
for  $\int \g(t)\dd W_H dt$ and $\int \g(t)\dd B_H(t)$ than for $\int \g(t)^2dt$, such as is
given in Proposition \ref{propWR}. 
It can be noted that combined  estimates from Proposition \ref{propWR} would
lead to estimate   $\E| I_H(\g)|\le \const \left(\E\int_{s}^{T}\g(t)^2dt\right)^{1/2}$. 
which is weaker than known estimates \citet{DecU,PT1}. \end{remark}
\begin{proposition} \label{propConst}
We have  that
\baaa
\int_{s}^{T}1\cdot \dd B_H(t)=B_H(T)-B_H(s).
\eaaa
 \end{proposition}
\subsubsection*{Extension on piecewise-predictable integrands from $\X_\e$}
\begin{definition}\label{intD} Let $\g\in\X_\e$, where $\e>0$.  By the definitions, there exists
a finite set $\TT$ of non-random times $\TT=\{T_k\}_{k=0}^n\subset [s,T]$, where $n>0$ is an integer, $T_0=0$, $T_n=T$,
and $T_{k+1}\in(T_k,T_k+\e]$ such that $\g(t)$ is $\GG_{T_k}$-measurable for $t\in[T_k,T_{k+1}]$. Let $\int_{T_{k-1}}^{T_k}\g(t)\dd B_H(t)$ be defined according to Definition  \ref{defWR} with the interval $[s,T]$ replaced by $[T_{k-1},T_k]$.
 We call the sum
\baaa
I_H(\g)=\sum_{k=1}^n\int_{T_{k-1}}^{T_k}\g(t)\dd B_H(t).
\eaaa
the foresighted   integral of $\g$ and denote it as $\int_{0}^{T}\g(t)\dd B_H(t)$.
\end{definition}
The integral in the above definition belongs  to $L_1(\O,\GG_T,\P)$ thanks to Propositions    \ref{propW} and \ref{propM}.
\begin{remark} \label{remPC}
It follows from Proposition \ref{propConst} that
\baaa
\int_{0}^{T}\g(t)\dd B_H(t)=\int_{0}^{T}\g(t)d_\RS B_H(t)
\eaaa
for piecewise constant  $\g\in\cup_{\e>0} \X_{\e,\PC}$.
However,  it appears that converges of Riemann sums  requires more restriction for non-piecewise constant $\g$ than the convergence
for the suggested new integral.. This is because
 this approximation  is finer that approximation by the piecewise constant functions.
\end{remark}

\subsection{Extension on random integrands of a general type with a mild restriction on prediction error}

Let $\E_t$ and $\Var_t$ denote the conditional expectation and the conditional variance given $\GG_t$, respectively

For $\nu> 0$ and $\e>0$, let  $\Y_{\nu,\e}$ be  the set of all processes $\g\in\L_{22}$ such that
 \baaa
&&\sup_{\tau\in[0,T]}\hspace{2mm}\sup_{t\in[\tau,T\land(\tau+\e) ]}\hspace{2mm}
[\E\Var_\tau \g(t)]^{1/2}\le C (t-\tau)^{1-H+\nu}\quad\hbox{a.s.}
\eaaa
for some $C=C(\g)>0$.

It can be noted that $\E_\tau\g(t)$ can be interpreted as the forecast at time $\tau$ of $\g(t)$ for $t>\tau$; the forecast is based on observations of the events from $\GG_\tau$. Respectively,
$\Var_\tau\g(t)$ can be interpreted as the conditional means-square error of this forecast given $\GG_\tau$.

In particular,  processes from  $\Y_{\nu,\e}$ with $\nu>0$  feature stronger predictability
on the short horizon $\e$  than processes from $\Y_{0,\e}$.

\begin{proposition}\label{propC} For any $\nu> 0$ and $\e>0$, the space $\Y_{\nu,\e}$ with the norm
\baaa
\|\g\|_{\Y_{\nu,\e}}\defi  \|\g\|_{\L_{22}}+
\sup_{\tau\in[0,T]}\hspace{2mm}\sup_{t\in[\tau,T\land(\tau+\e) ]}\hspace{2mm}
\left[\E\Var_\tau \g(t)\right]^{1/2}/(t-\tau)^{1-H+\nu}.
\eaaa
is a Banach space.
\end{proposition}

It follows from the definitions that if $\e_0\in (0,\e)$ and $\g\in\Y_{\nu,\e}$ then
 $\g\in\Y_{\nu,\e_0}$ and $\|\g\|_{\Y_{\nu,\e_0}}\le \|\g\|_{\Y_{\nu,\e}}$.
 Also, it can be seen that $\X_\e\subset \Y_{\nu,\e}$ for any $\nu>0$.

Let  $\Y\defi \cup_{\nu>0, \e>0}\Y_{\nu,\e}$.

Clearly, the set $\Y$ is everywhere dense in $\L_{22}$.

\begin{example} We have that $B|_{[0,T]}\in \Y_{0,\e}$ but $B|_{[0,T]}\notin \Y$.
On the other hand, $B_H|_{[0,T]}\in \Y_{{2H-1},\e}$ for any $\e>0$.
\end{example}

For  $\g\in\L_{22}$,  let $\Z(\g)$ be the set of processes $\{\g_n\in\X_{},\quad n=0,1,2,...\}$, such that $\g_n(t)=\E_{T_k}\g(t)$ for $t\in [T_k,T_{k+1})$, where $ k=0,1,...,2^n$ and where $T_{k}=kT/2^{n}$.

\begin{theorem}\label{ThM}  \begin{enumerate}
\item Let   $\g\in\Y_{}$, and let
$\{\g_n\}_{n=1}^\infty=\Z(\g)$. Then the sequence $\{I_H(\g_n)\}_{n=1}^\infty$ converges  to a limit
 in $ L_1(\O,\GG_T,\P)$ uniformly over
 $H\in(1/2, c)$ for any $c\in (1/2,1)$. Let  $I_H(\g)$ denote this limit. 
\item For any $\e>0$, $H\in (1/2,1)$, and $\nu>0$,  the  operator $I_H(\cdot): \Y_{\nu,\e}  \to L_1(\O,\GG_T,\P)$ defined in statement (i) is a linear continuous operator.
For any $\e>0$, the norms of these operators are bounded  in $H\in (1/2,c)$, for any $c\in (1/2,1)$.
\end{enumerate}
\end{theorem}
We will regard $I_H(\g)$ defined in Theorem \ref{ThM}  as the stochastic integral
\baaa
I_H(\g)=\int_0^T\g(t)\dd B_H(t),\quad \g\in\Y_0.
\eaaa
\begin{corollary}\label{corrI}
For any $\e>0$ and $\nu> 0$, there exists a constant $c>0$ depending on $T$,$\e$,$\nu$ only
such that
\baaa
\E\left|\int_0^T\g(t)\dd B_H(t)\right|\le c\|\g\|_{\Y_{\nu,\e}}\quad \forall \g\in \Y_{\nu,\e}.
\eaaa
\end{corollary}
Corollary \ref{corrI} follows immediately  from Theorem \ref{ThM}.

For $\nu>0$ and $r>1$, let  $\H_{\nu,r}$ be the set of all $\g\in \L_{22}$  such that
 $\sup_{s,t\in[0,T]}\|\g(s)-\g(t)\|_{L_r(\O)}\le C|t-s|^{1-H+\nu}$ for some $C=C(\g)>0$.

 It can be seen that $\H_{\nu,r}\subset \Y_{\nu,\e}$ for $r\ge 2$ for all $\e>0$.

 For  $\g\in\H_{\nu.r}$,  let $\oo\Z(\g)$ be the set of processes
$\{\g_n\in\X_{},\quad n=0,1,2,...\}$, such that,   for $t\in [T_k,T_{k+1})$,   either
$\g_n(t)=\g(T_k)$, or $\g_n(t)=\E_{T_k}\g(t)$,  where $ k=0,1,...,2^n$ and where $T_{k}=kT/2^{n}$.

\begin{proposition}\label{ThHolder} For any $r\in (1,2]$ and  $\nu>0$,
the conclusions of Theorem \ref{ThM} hold for $\g\in\H_{\nu,r}$ if  $\Y_{}$, 
$\Y_{\nu,\e}$, and $\Z(\g)$,    are replaced by $\cup_{\nu>0}\H_{\nu,r}$, $\H_{\nu,r}$,  and $\oo\Z(\g)$,  respectively.
\end{proposition}

  \section{Continuity of the foresighted  integral in $H\to 1/2+0$}
  \label{SecC}
The following theorem describes some classes of random integrands where
the stochastic integrals are continuous  with respect to
the Hurst parameter $H\to 1/2+0$.
 \begin{theorem}\label{ThL}
For any $\g\in\Y$,  \baa
\E\left|\int_0^T\g(t)\dd B_H(t)-\int_0^T\g(t)dB(t)\right|\to 0 \quad\hbox{as}\quad H\to 1/2+0.
\label{lim}\eaa
\end{theorem}

In fact, the question about continuity at  $H\to 1/2$ of  stochastic integrals with respect to $dB_H$ is quite interesting.
In particular, it is known
 that
\baa
\E\int_0^TB_H(t)d_{\RS}B_H(t)\nrightarrow \E\int_0^TB(t)dB(t) \quad\hbox{as}\quad H\to 1/2+0.
\label{nonlim}\eaa
This follows from the equality \baaa
2\int_0^T B(t)dB(t)=B(T)^2-T
\label{sqBB}\eaaa
combined with the  equalities  \cite{Shi}
\baaa
2\int_0^TB_H(t)d_{\RS}B_H(t)=B_H(T)^2,\quad  H\in(1/2,1).
\label{sqB}\eaaa
  \begin{remark}\label{remNonprot}
Theorem \ref{ThL} does not contradict to the divergence stated  in  (\ref{nonlim})
since $B_{[0,T]}\notin \Y$. On the other hand, this theorem ensures that, for any $H_1>1/2$,
\baaa
\E\int_0^TB_{H_1}(t)\dd B_H(t)\rightarrow \E\int_0^TB_{H_1}(t)dB(t) \quad\hbox{as}\quad H\to 1/2+0,
\eaaa
 since $B_H|_{[0,T]}\notin\Y_{}$.
\end{remark}

 \section{Proofs}\label{SecP}
 Consider the derivative
\baaa
f_t'(t,r)= \ddH (t-r)^{H-3/2},\quad t>r.
\eaaa
Since $H-3/2\in (-1,-1/2)$, it follows that  $2(H-3/2)\in (-2,-1)$ and  $\|f'_{t}(t,\cdot)\|_{L_2(-\infty,\Ts)}<+\infty$ for all  $\Ts<t$.

{\em Proof of Proposition \ref{propW}}.
For $\tau\in [s,T]$, $\e\ge 0$,  and  $g\in L_2(s,T)$, set \baa
G_{H,\e}(\tau,T,g)\defi  \ddH \int_{\tau}^T(t-\tau+\e)^{H-3/2}g(t)dt.
\eaa

By the restrictions on $\g$ and by (\ref{estG}), we have that $G_H(\cdot,T,\g)$
is $\GG_s$-measurable for any $\tau$, that
$\int_{s}^T dB(\tau)G_{H}(\tau,T,\g)$ is well defined as an It\^o's integral, and that
$\int_{s}^T \g(t)dW_{H,\e}(\tau)$ is also well defined as the It\^o's integral
\baa
&&\int_{s}^T \g(t)dW_{H,\e}(t)\nonumber\\&&=c_H\int_{s}^T \g(t)f(t,t-\e)dB(t) + \ddH \int_{s}^T \g(t)dt\int_s^t (t-\tau+\e)^{H-3/2}dB(\tau)\nonumber\\&&= \ddH \int_{s}^T dB(\tau)\int_\tau^T(t-\tau)^{H-3/2}\g(t)dt,
\label{IW0}\eaa
i.e.
\baa
\int_{s}^T \g(t)dW_{H,\e}(t)=\int_{s}^T dB(\tau)G_{H,\e}(\tau,T,\g).
\label{IW}\eaa
Furthermore, let
\baaa D_\e\defi \int_{s}^T dB(\tau)G_{H}(\tau,T,\g)-\int_{s}^T \g(t)dW_{H,\e}(t).
 \eaaa
We have that $D_\e=\oo D_\e +\w D_\e$, where
$\oo D_\e\defi \int_{s}^T \g(t)f(t,t-\e)dB(t)$ and  where
\baaa \w D_\e
\defi \int_{s}^T dB(\tau)[G_{H}(\tau,T,\g)-G_{H,\e}(\tau,T,\g)].
 \eaaa
Clearly, $\E\oo D_\e^2\to 0$ as $\e\to 0$. Let us show that  $\E\w D_\e^2\to 0$ as $\e\to 0$.

It suffices to consider $\e=\e_j$ for a monotonically decreasing sequence $\{\e_j\}_{j=1}^{\infty}$.

Assume first that  $\g(t)\ge 0$ a.e..
In this case, $(t-\tau+\e_i)^{H-3/2}\g(t)> (t-\tau+\e_j)^{H-3/2}\g(t)\ge 0$ a.e. if  $i>j$, i.e., $\e_i<\e_j$.

It follows that $G_{H}(\tau,T,\g)-G_{H,\e}(\tau,T,\g)\ge 0$ a.s. for almost all
$\tau$. It also follows $\|G_{H,\e}(\cdot,T,\g)\|_{L_2(s,T)}\le c\|\g\|_{L_2(s,T)}$
with the same $c$ as in (\ref{estG}).

We have  that $G_{H}(\tau,T,\g)-G_{H,\e}(\tau,T,\g)\to 0$  a.s. for almost all
$\tau$ as $\e=\e_j\to 0$ and that
$0\le G_{H,\e}(\tau,T,\g)\le G_{H}(\tau,T,\g)$ for a.e. $\tau$.  By the Lebesgue
Dominated Convergence Theorem, it follows that $\E\w D_\e^2\to 0$ as $\e\to 0$.

The case where $\g\le 0$ can be considered similarly. In the case of a sign variable $\g$, apply the proof above for $\g_+=\g\Ind_{\g\ge 0}$ and for $\g_-=-\g\Ind_{\g\le  0}$ separately. Then the proof for $\g=\g_+-\g_-$ follows.
This completes the proof of Proposition \ref{propW}.

{\em Proof of Proposition \ref{propM}}.
Let us prove statement (i).  We need to verify the properties related to the differentiability of
 $R_H(t)$.

Let $t>\Ts$ and $r<s$.

Let $f^{(1)}(t,r,\d)\defi (f(t+\d,r)-f(t,r))/\d$, where $\d\in (-(t-\Ts)/2,(t-\Ts)/2)$.

Clearly, $f'_t(t,r)-f^{(1)}(t,r,\d)\to 0$ as $\d\to 0$ for all $t>s$ and $r<s$.
Let us show that $\|f'_t(t,\cdot)-f^{(1)}(t,\cdot,\d)\|_{L_2(-\infty,\Ts)}\to 0$ as $\d\to 0$.
We have that
\baaa
f^{(1)}(t,r,\d)=\d^{-1}\int_t^{t+\d}f_t'(s,r)ds=f'_t(\t(t,\d),r)
\eaaa
for some $\t(t,\d)\in (t,t+\d)$.  Hence
\baa
|f'_t(t,r)-f^{(1)}(t,r,\d)|\le \sup_{h\in (t,t+\d)}|f'_t(t,r)-f'_t(h,r)|\le  \d\sup_{h\in (t,t+\d)}|f''_{tt}(h,r)|,
\label{est}\eaa
where
\baaa
f''_{tt}(h,r)= \ddH (H-3/2)(h-r)^{H-5/2}.
\eaaa
For $\d>0$, we have that
\baaa
\sup_{h\in (t,t+\d)}|f''_{tt}(h,r)|\le d_H|(H-3/2)|(t-r)^{H-5/2}.
\eaaa
For $\d\in (-(t-\Ts)/2,0]$, we have that
\baaa
\sup_{h\in (t,t+\d)}|f''_{tt}(h,r)|\le d_H|(H-3/2)|(t+\d-r)^{H-5/2}.
\eaaa
It follows that $\|f''_{tt}(t,\cdot)\|_{L_2(-\infty,\Ts)}<+\infty$.

By (\ref{est}), it follows for all $t>s$ \baaa
\DD R_H(t)=\lim_{\d\to 0} \frac{R_H(t+\d)-R_H(t)}{\d}=\cH   \int_{-\infty}^{\Ts}f'_t(t,r)dB(r),
\eaaa
for  the mean square limit described in  statement (ii).

Further, we have that
\baa
\E \DD R_H(t)^2&=&\cH \int_{-\infty}^{\Ts}|f_t'(t,r)|^2dr
=d_H^2\int_{-\infty}^{\Ts}(t-r)^{2H-3}dr\nonumber\\&=&\frac{d_H^2}{2-2H}(t-r)^{2H-2}\Bigl|^{\Ts}_{-\infty}
=\frac{d_H^2}{2-2H}(t-\Ts)^{2H-2}.
\label{DRn}
\eaa
Hence, for $t>\Ts$,
\baaa
\E\int_{\Ts}^t \DD R_H(r)^2dr =\frac{d_H^2}{2-2H}\int_{\Ts}^t (r-\Ts)^{2H-2} dr&=&
\frac{d_H^2}{(2-2H)(2H-1)}(t-\Ts)^{2H-1}\nonumber\\&=& \frac{c_Hd_H}{2(2-2H)}(t-\Ts)^{2H-1}.
\label{EDR1} \eaaa
This completes the proof of Proposition \ref{propM}.  $\Box$

{\em Proof of  Proposition \ref{propWR}} follows from (\ref{estG}) and Proposition \ref{propM}. $\Box$

{\em Proof or Proposition \ref{propConst}}. Let \baaa
h\defi H-1/2. 
\eaaa
By the definitions,
\baaa
&&\int_{s}^{T}1\cdot \dd B_H(t)
=J_1+J_2,
\eaaa
where
\baaa
J_1\defi \int_{s}^{T}1\cdot d_F W_H(t)=   d_H\int_{s}^{T}d B(\tau) \int_\tau^T (t-\tau)^{h-1}dt =\int_{s}^{T}d B(\tau)G_H(\tau,T,1)
\eaaa
and
\baaa
J_2\defi \int_{s}^{T}1\cdot \DD R_H(t)dt=\int_s^T dt \int_{-\infty}^s f'_t(t,r) dB(r)=  d_H\int_s^T dt \int_{-\infty}^s (t-\tau)^{h-1}dB(\tau).
\eaaa
We have that
\baaa
J_1=c_H\int_{s}^T dB(\tau)  (T-\tau)^{h}
\eaaa
and
\baaa
J_2=d_H\int_{-\infty}^s dB(\tau) \int_{s}^T (t-\tau)^{h-1}dt=c_H\int_{-\infty}^s dB(\tau)  [(T-\tau)^{h}-(s-\tau)^h].
\eaaa
Hence
\baaa
\int_{s}^{T}1\cdot \dd B_H(t)=J_1+J_2=c_H\int_{s}^T dB(\tau)  (T-\tau)^{h}
+c_H\int_{-\infty}^s dB(\tau)  [(T-\tau)^{h}-(s-\tau)^h].
\eaaa
It follows  from the well known properties of fractional Brownian motions that this value is
$B_H(T)-B_H(s)$. Let us show this for the sake of completeness. We have that
\baaa
&&B_H(T)-B_H(s)=c_H\int_{0}^T dB(\tau)  (T-\tau)^{h}+c_H\int_{-\infty}^0 dB(\tau)  [(T-\tau)^{h}-(-\tau)^h]\\&&
\hphantom{xxxxxxx}-c_H\int_{0}^s dB(\tau)  (s-\tau)^{h}-c_H\int_{-\infty}^0 dB(\tau)  [(s-\tau)^{h}-(-\tau)^h]
\\
&&
=c_H\int_{s}^T dB(\tau)  (T-\tau)^{h}
+c_H\int_{-\infty}^s dB(\tau)  [(T-\tau)^{h}-(s-\tau)^h].\eaaa
This completes the proof of Proposition \ref{propConst}. $\Box$

 {\em Proof of Proposition \ref{propC}}.
We denote by $\oo\ell _{1}$ the Lebesgue measure in $\R$, and we
denote by $ \oo{{\cal B}}_{1}$ the $\sigma$-algebra of Lebesgue
sets in $\R$. Let $D=\{(t,r):0\le r\le t\le T\}$.
\par
 Let $\V_1=L_2([0,T], \oo{{\cal B}}_{1},\oo\ell _{1}, L_2(\O,\GG_0,\P))$, and
let $\V_2$ be the linear normed
 space of all
measurable function (classes of equivalency) $g:D\times\O\to\R$ such that
$g(t,r)\in L_2(\O,\GG_r,\P)$ for a.e. $t,r$,
with the norm
 \baaa
\|\w g\|_{\V_2}&=&\left(\E\int_0^Tdt\int_0^t g(t,r)^2 dr\right)^{1/2} \\&+&\sup_{\tau\in[0,T]}\hspace{2mm}
 \sup_{t\in [\tau,(\tau+\e)\land T]}\left(\E \int_\tau^t g(t,r)^2 d\t\right)^{1/2}/(t-\tau)^{1-H+\nu}.\eaaa
By Clark's theorem, it follows that $\g\in\Y_{\nu,\e}$  can be represented as
\baaa
\g(t)=\E_0\g(t)+\int_{0}^tg(t,r)dB(r)
\eaaa
for some $g(t,r)\in \V_2$; here $\E_0\g(t)\in \V_1$.  In this case,  $\Var_{\tau} \g(t)=\E_\tau\int_\tau^t g(t,r)^2dr$.
To prove the proposition, it suffices to observe that the space
$\V_1\times \V_2$ is complete and  is in a continuous and continuously invertible bijection with the space $\Y_{\nu,\e}$. This completes the proof of Proposition \ref{propC}. $\Box$
\par
\vspace{2mm}
\par
To prove  Theorems \ref{ThM}, Proposition \ref{ThHolder}, and Theorem \ref{ThL},  we will need some notation.

We will be using functions
\baa
\w \rho(t)\defi \cH \int_{-\infty}^0 f'_t(t,r)dB(r), \quad \rho(t,\tau)\defi \cH \int^{\tau}_0 f'_t(t,r)dB(r), \quad \tau> t>0.
\label{rho}
\eaa

In the proofs below, we consider an integer $n>0$ and $\g_n\in\X_{}$ such  that
there exist some $\e>0$ and a set $\TT_{\g_n}=\{T_k\}_{k=1}^n\subset[0,T]$,   where $T_0=0$, $T_n=T$,
and $T_{k+1}\in(T_k,T_k+\e)$ such that $\g_n(t)\in L_2(\O,\GG_{T_k},\P)$ for $t\in[T_k,T_{k+1})$.
\par
Let
\baaa
I_{W,H,k}=\int_{T_{k-1}}^{T_k}\g_n(t)dW_{H,k}(t),\quad  I_{R,H,k}=\int_{T_{k-1}}^{T_k}\g_n(t)\DD R_{H,k}(t)dt,
\eaaa
 where $W_{H,k}$, $R_{H,k}$, and $\DD R_{H,k}$ are defined similarly to $W_H$, $R_H$, and
$\DD R_H$, with $[s,T]$ replaced by $[T_{k-1},T_k]$.

Let  \baa
I_{W,H}(\g_n)\defi \sum_{k=1}^{n} I_{W,H,k},\quad I_{R,H}(\g_n)\defi \sum_{k=1}^{n} I_{R,H,k}.
\label{IHGW}\eaa

Clearly,
\baaa
I_H(\g_n\Ind_{[T_{k-1},T_k)})=\int_{T_{k-1}}^{T_k}\g_n(t)\dd B_H(t)=I_{W,H,k}+I_{R,H,k},
\eaaa
and
\baaa
I_H(\g_n)=I_{W,H}(\g_n)+I_{R,H}(\g_n).
\label{IHG}\eaaa
\par

By the definitions,
\baaa
I_{R,H,k}&= &\int_{T_{k-1}}^{T_{k}}\g_n(t)\DD R_k(t)dt
=\cH \int_{T_{k-1}}^{T_{k}}\g_n(t)\int_{-\infty}^{T_{k-1}}f'_t(t,s)dB(s)
\\&= &\cH \int_{T_{k-1}}^{T_{k}}\g_n(t)\int_{-\infty}^{0}f'_t(t,s)dB(s) + \cH \int_{T_{k-1}}^{T_{k}}\g_n(t)\int_{0}^{T_{k-1}}f'_t(t,s)dB(s)
\\&= &\int_{T_{k-1}}^{T_{k}}\g_n(t)\w \rho(t)dt+\int_{T_{k-1}}^{T_{k}}\g_n(t)\rho(t,T_{k-1})dt.
\eaaa
Hence
\baaa
I_H(\g_n)=I_{W,H}(\g_n)+\w I_{R,H}(\g_n)+\oo J_{R,H}(\g_n),
\eaaa
where
\baa
\w I_{R,H}(\g_n)=\int_0^T\g_n(t)\w  \rho(t)dt, \qquad
\oo J_{R,H}(\g_n)\defi  \sum_{k=1}^n J_{R,H,k},  \label{IJ}
\eaa
\baa
J_{R,H,k}=\int_{T_k}^{T_{k+1}}\g_n(t)\rho(t,T_k)dt.\label{IIJ}
\eaa
\par
For $k=0,...,n-1$,
consider operators $\G_k(\cdot):L_2(0,T_{k+1})\to L_2(0,T_{k+1})$ such that
\baaa
\G_k(\cdot, g)=G_H (\cdot,T_{k+1},g),\label{GGH}
\label{GGHexpl}
\eaaa
i.e.
\baa
\G_k(\tau, g)= d_H \int_{\tau}^{T_{k+1}}(t-\tau)^{H-3/2}g(t)dt.
\label{GH}\eaa

By the properties of the Riemann--Liouville integral, $
\|\G_k(\cdot,g)\|_{L_2(T_k,T_{k+1})}\le\w c\|g\|_{L_2(T_k,T_{k+1})}$   for some $\w c>0$ that
is independent on $g\in L_2(T_k,T_{k+1})$ and  $H\in(1/2,1)$.
\par
\begin{lemma}\label{lemmaW}  For any $c\in(1/2, 1)$, there exists some $C=C(c)>0$  such that, for any $\g_n\in\X_{}$ and $H\in (1/2,1)$,
\baa
\E|I_{W,H}(\g_n)|+\E|\w I_{R,H}(\g_n)|\le C \|\g_n\|_{\L_{22}}. \eaa
\end{lemma}
\par
{\em Proof of Lemma \ref{lemmaW}}.
For $k=1,...,n$, we have that
\baaa
&&I_{W,H,k}= \ddH \int_{T_{k-1}}^{T_k} \g_n(t)dt\int_{T_{k-1}}^t (t-\tau)^{H-3/2}dB(\tau)\\&&= \ddH \int_{T_{k-1}}^{T_k}  dB(\tau)\int_\tau^{T_k}(t-\tau)^{H-3/2}\g_n(t)dt
=\int_{T_{k-1}}^{T_k}  dB(\tau)\G_{k-1}(\tau,\g_n).
\eaaa
The last  integral here converges  in $L_2(\Omega,\GG_{T},\P)$.
Hence
\baaa
\E\|I_{W,H}(\g_n)\|_{L_2(\O)}^2=\E \left(\sum_{k=1}^n I_{W,H,k}\right)^2=\sum_{k=1}^n\E I_{W,H,k}^2=
\E\sum_{k=1}^n\int_{T_{k-1}}^{T_k} \G_{k-1}(\tau,\g_n)^2d\tau
\\\le \w c\,\E\sum_{k=1}^n\int_{T_{k-1}}^{T_k} \g_n(\tau)^2d\tau= \w c\, \|\g_n\|^2_{\L_{22}}.
\eaaa
\par
Further, we have that
\baaa
\E|\w I_{R,H}(\g_n)|\le \left(\E\int_{0}^{T}\g_n(t)^2dt\right)^{1/2} \left(\E\int_{0}^{T}\w\rho(t)^2dt\right)^{1/2}.
\eaaa
By (\ref{EDR}),
$\E\int_{0}^{T}\w\rho(t)^2dt\le \frac{d_H^2}{2(2-2H)}T^{2H-1}$.
This completes the proof of Lemma \ref{lemmaW}. $\Box$

The following proofs will be given for Theorem \ref{ThM} and  \ref{ThL} simultaneously
with the proof of Proposition \ref{ThHolder}.  

For the sake of the proofs of Theorem \ref{ThM} and  \ref{ThL},
we assume below that that $r=2$, $p=2$,
$\g\in \Y_{\nu,\e}$  and $\{\g_n\}_{n=1}^\infty=\Z(\g)$.
For the sake of the proof of Proposition \ref{ThHolder}, we assume below that $r\in (1,2]$, $p=(1-1/r)^{-1}$,
$\g\in \H_{\nu,r}$ and  $\{g_n\}_{n=1}^\infty =\oo\Z(\g)$.

We consider below positive integers $n,m\to +\infty$ such that $n\ge m$. We assume below that $T_k=kT/2^n$, $k=0,1,...,2^n$. This means that   the grid $\{T_k\}_{k=0}^{2^n}$
is formed as defined for $n$ rather than for $m$; since $n\ge m$, Definition \ref{intD} is applicable to the integral $\int_0^T\g_m(t)\dd B_H(t)$
with this grid as well.

We denote 
\baaa
\e_m\defi T/2^n=T_{k	+1}-T_k, \e_n\defi T/2^n=T_{k	+1}-T_k.
\eaaa

We assume that $m$ is such that $\e_m\le \e$. It implies that $\e_n\le \e$ as well.

We denote by   $J_{R,k,m}$ and $J_{R,k,n}$ the corresponding values  $J_{R,H,k}$ defined for $\g=\g_n$  and $\g=\g_m$ respectively obtained using  the
same grid  $\{T_k\}_{k=0}^{2^n}$.

\begin{lemma}\label{lemmaR} The
sequence $\{I_{R,H}(\g_n)\}_{n=1}^\infty$ has a limit in
$L_1(\O,\GG_T,\P)$; it converges to this limit uniformly in $H\in (1/2,c)$, for any $c\in (1/2,1)$.
\end{lemma}
\par
{\em Proof of Lemma \ref{lemmaR}}.
Clearly,
\baa
\|\g_n-\g_m\|^2_{\L_{22}}\to 0\quad\hbox{as}\quad n,m\to +\infty
\label{g_nm}
\eaa
and 
\baa
\|\g_n-\g_m\|^2_{\L_{22}}\to 0\quad\hbox{as}\quad m\to +\infty \quad
\hbox{uniformly in }\ n>m.
\label{g_mnu}
\eaa

By Lemma \ref{lemmaW}, we have that
\baaa
\E\|I_{W,H}(\g_n)-I_{W,H}(\g_m)\|_{L_2(\O)}^2+
\E|\w I_{R,H}(\g_n)-\w I_{R,H}(\g_m)|\to 0\quad\hbox{as}\quad b,m\to +\infty.
\eaaa
This implies that the sequences
$\{I_{W,H}(\g_n)\}_{n=1}^\infty$ and $\{\w I_{R,H}(\g_n)\}_{n=1}^\infty$ have limits in
$L_1(\O,\GG_T,\P)$, and that  they converge to these limits uniformly in $H\in (1/2,c)$, for any $c\in (1/2,1)$.

Therefore, to prove Lemma \ref{lemmaR}, it suffices  to prove that the sequence
$\{\oo J_{R,H}(\g_n)\}_{n=1}^\infty$ have a limit in
$L_1(\O,\GG_T,\P)$ as well, and that  it converges to this limit uniformly in $H\in (1/2,c)$, for any $c\in (1/2,1)$.
.

Let \baaa
\xi_k(t)\defi \rho(t,T_k)=d_H\int_0^{T_k}(t-s)^{H-3/2}dB(s).
\label{kk}
\eaaa

  We have that
\baaa
\psi_{n,m,k}\defi J_{R,k,n} -J_{R,k,m}=  \int_{T_k}^{T_{k+1}}[\g_n(t)-\g_m(t)]\xi_k(t)dt,
\label{kk0}
\eaaa

Remind that $p>0$ is such that $1/p+1/r=1$. We have that
\baaa
\|\psi_{n,m,k}\|_{L_1(\O)} \le \int_{T_k}^{T_{k+1}}\|\g_n(t)-\g_m(t)\|_{L_r(\O)}\|\xi_k(t)\|_{L_p(\O)}dt.\label{kk2}
\eaaa
Further, we have that
\baaa
\|\xi_k(t)\|_{L_2(\O)}^2=d_H^2\int_{0}^{T_{k}}(t-s)^{2H-3}ds
=\frac{d_H^2}{2H-2}\left[(t-T_k)^{2H-2}-t^{2H-2}\right]\\
=\frac{d_H^2}{2-2H}\left[t^{2H-2}-(t-T_k)^{2H-2}\right], \quad t\in (T_k,T_{k+1}].
\label{kk3}\eaaa
Hence
\baaa
\int_{T_k}^{T_{k+1}}\|\xi_k(t)\|_{L_2(\O)}^2 dt&=&
\frac{d_H^2}{(2-2H)(2H-1)} [(T_{k+1}-T_k)^{2H-1}-T_{k+1}^{2H-1}+T_k^{2H-1}]
\\&=&\frac{c_H d_H }{4-4H}[(T_{k+1}-T_k)^{2H-1}-T_{k+1}^{2H-1}+T_k^{2H-1}].
\label{kk33}
\eaaa
 Hence
\baa
&&\left(\int_{T_k}^{T_{k+1}}\|\xi_k(t)\|_{L_2(\O)}^2 dt\right)^{1/2}\le \oo C_0 C_H\e_n^{H-1/2},
\label{kk3333}
\eaa
where
\baa
C_H\defi \frac{\sqrt{c_H d_H }}{2-2H},
\label{CH}
\eaa   and where
$\oo C_0>0$ is independent on $\g$, $k$ and $H$; it depends on $T$ only.

By the properties of Gaussian distributions, we have that
\baaa
\|\xi_k(t)\|_{L_p(\O)}\le  C(p)\|\xi_k(t)\|_{L_2(\O)}\eaaa
for some $C(p)>0$. Hence
\baa
&&\int_{T_k}^{T_{k+1}}\|\xi_k(t)\|_{L_p(\O)}dt
\le C(p) \int_{T_k}^{T_{k+1}}\|\xi_k(t)\|_{L_2(\O)}dt\le C(p) \left(\int_{T_k}^{T_{k+1}}\|\xi_k(t)\|_{L_2(\O)}^2dt \right)^{1/2}\e_n^{1/2}
\nonumber\\&&
\le C(p)\oo C_0  C_H\e_n^{H-1/2}\e_n^{1/2}=C(p)\oo C_0  C_H\e_n^{H}.
\label{kk33}\eaa
\par
Let
\baaa
 T^{(m)}_d\defi \e_m d, \quad d=0,1,...,2^m,
\label{em}\eaaa
Here $\e_m= T/2^{m}$. Let 
\baaa
\tau_m(t)\defi \inf\{T^{(m)}_d:\ t\in[T^{(m)}_d,T^{(m)}_{d+1}),\ d=0,1,...,2^m-1\}.\index{,\qquad
 \oo\tau_m(t)=\inf\{T^{(m)}_d:\ t\in(T^{(m)}_{d-1},T^{(m)}_{d}]\}.}
\eaaa
Clearly,  the function $\tau_m(t)$ is non-decreasing, and  $\tau_m(t)\le \tau_n(t)$.
\par
By the definitions, we have that $\g_m(t)=\E_{\tau_m(t)}\g(t)
=\E_{\tau_m(t)}\g_n(t)$ and $\g_n(t)=\E_{\tau_n(t)}\g(t)$. Hence 
\baaa
\|\g_n(t)-\g_m(t)\|_{L_2(\O)}=\|\g_n(t)-\E_{\tau_m(t)}\g_n(t)\|_{L_2(\O)}\le
\|\g(t)-\E_{\tau_m(t)}\g(t)\|_{L_2(\O)}.
\eaaa
For the sake  of the proof of Theorem \ref{ThM},  we have assumed that $\g\in \Y_{\nu,\e}$. It follows that 
\baa
\sup_{t\in [0,T]}\|\g_m(t)-\g_n(t)\|_{L_2(\O)}
\le
 \sup_{t\in [0,T]} \left(\E\Var_{\tau_m(t)}\g(t)\right)^{1/2} \le c\e_m^{1-H+\nu}\|\g_m\|_{\Y_{\nu,\e}},
 \label{gg}
\eaa
where $c>0$ are independent on $\g$ and $H\in (1/2,1)$.

\par
Let $n=m+1$. In this case,  we have that $\e_m=2\e_n$. By (\ref{kk3333}) and (\ref{gg}), we have that and
\baaa
&&\|\psi_{k,m+1,m}\|_{L_1(\O)}\le \int_{T_k}^{T_{k+1}}\|\g_{m+1}(t)-\g_m(t)\|_{L_r(\O)}
\|\xi_k(t)\|_{L_p(\O)}dt
\\&&\le 
\sup_{t\in[0,T]}\|\g_{m+1}(t)-\g_m(t)\|_{L_r(\O)}
\int_{T_k}^{T_{k+1}}
\|\xi_k(t)\|_{L_p(\O)}dt
\\&&\le \int_{T_k}^{T_{k+1}}\|\g_{m+1}(t)-\g_m(t)\|_{L_r(\O)}
\|\xi_k(t)\|_{L_p(\O)}dt
\\&&\le c_\psi C_H\e_m^{1+\nu}\|\g_m\|_{\Y_{\nu,\e}},
\label{bb1}
\eaaa
where $c_\psi>0$ is independent on $\g$, $k$, and $H\in (1/2,1)$. We have that $2^n=2^{m+1}=2 T/\e_m$. Hence
\baa
\|\oo J_{R,H}(\g_{m+1})-\oo J_{R,H}(\g_m)\|_{L_1(\O)}
\le \E\sum_{k=0}^{2^n-1}\|\psi_{k,n,m}\|_{L_1(\O)}
\le 2^n c_\psi C_H  \e_m^{1+\nu} \|\g\|_{\Y_{\nu,\e}}\nonumber\\=2 T\e_m^{-1} c_\psi C_H \e_m^{1+\nu} \|\g\|_{\Y_{\nu,\e}}
  =c_J C_H
(2^{-m})^{\nu}\|\g\|_{\Y_{\nu,\e}},
\label{J1}\eaa
where   $c_J>0$ is independent on $m$, $\g$, $H$, and $\nu\ge 0$.

Further, let $m\in\{1,2,...,n\}$.   We have that
\baa
&&\oo J_{R,H}(\g_n)-\oo J_{R,H}(\g_m)\nonumber\\&=&\oo J_{R,H}(\g_n)-\oo J_{R,H}(\g_{n-1})+\oo J_{R,H}(\g_{n-1})-\oo J_{R,H}(\g_m)\nonumber\\&=&
\oo J_{R,H}(\g_n)-\oo J_{R,H}(\g_{n-1})+\oo J_{R,H}(\g_{n-1})-\oo J_{R,H}(\g_{n-2})+\oo J_{R,H}(\g_{n-2})-\oo J_{R,H}(\g_m)\nonumber\\&
=&...= \sum_{k=m+1}^n (\oo J_{R,H}(\g_{k})-\oo J_{R,H}(\g_{k-1})).
\label{J2}\eaa
It follows that
\baa
\|\oo J_{R,H}(\g_n)-\oo J_{R,H}(\g_m)\|_{L_1(\O)}\le c_JC_H\sum_{k=m+1}^{n} (2^{-k})^{\nu}\|\g\|_{\Y_{\nu,\e}}\to 0\quad \nonumber\\\hbox{as}\quad m\to +\infty
\label{J3}
\eaa
uniformly in $n>m$ and in the case where $\nu>0$, uniformly in $H\in (1/2,c)$, for any $c\in (1/2,1)$.
Hence $\{\oo J_{R,H}(\g_n)\}$ is a Cauchy sequence in $L_q(\O,\F,\P)$, and has a limit in this space, uniformly in $H\in (1/2,c)$, for any $c\in (1/2,1)$.

For the sake of the proof of Proposition \ref{ThHolder},
 we use, instead of (\ref{gg}), the estimates  
\baaa
&&\sup_t\|\g_n(t)-\g_m(t)\|_{L_r(\O)}
 =\sup_{k\in\{0,...,2^n-1\}} \sup_{t\in [T_k,T_{k+1}]}\|\g_m (t)-\g_n(t)\|_{L_r(\O)}\\
 &&\le  \e_m^{1-H+\nu}\|\g\|_{\H_{\nu,r}}.
\eaaa
Then the proof above can be repeated with minor changes. In particular, the corresponding constant $c_J$ depends on $r$.

This completes the proof of Lemma \ref{lemmaR}. $\Box$

{\em Proof of Theorem \ref{ThM}}. It
follows immediately from Lemma \ref{lemmaR} that the sequence
$\{I_{H}(\g_n)\}_{n=1}^\infty$ converges to a limit in
$L_1(\O,\GG_T,\P)$, uniformly in $H\in (1/2,c)$, for any $c\in (1/2,1)$.
This proves statement (i)  of Theorem \ref{ThM}.

Let us prove statement (ii) of Theorem \ref{ThM}.
It follows from Lemma  \ref{lemmaW} that the operators
$I_{W,H}(\cdot): \X_{}  \to L_1(\O,\GG_T,\P)$  and
$\w I_{R,H}(\cdot): \X_{}\to L_1(\O,\GG_T,\P)$ allow continuous extension
into continuous operators $I_{W,H}(\cdot): \L_{22}\to L_1(\O,\GG_T,\P)$  and
$\w I_{R,H}(\cdot): \L_{22}\to L_1(\O,\GG_T,\P)$, that are bonded uniformly in $H\in (1/2,c)$, for any $c\in (1/2,1)$.

It suffices to show that, for any $\nu>0$ and $\e>0$,  
\baaa
\sup_{n\ge 0}\E|\oo J_{R,H}(\g_{n})|\le \oo C \|\g\|_{\Y_{\e,\nu}}
\eaaa
for some $\oo C=\oo C(\e,\nu)>0$.

\par
Assume that $\g\in\Y_{\nu,\e}$ for some $\e>0$. Let
\baa
m_\e\defi\min\{m:\  2^{-m}T\le\e\}.
\label{m_e}
\eaa
 It follows from (\ref{J2}) that, for all  $n>m_\e$,
\baa
\E|\oo J_{R,H}(\g_n)|\le \|\oo J_{R,H}(\g_{m_\e})\|_{L_1(\O)}
+ c_JC_H\sum_{k=m+1}^{n} (2^{-k})^{\nu}\|\g\|_{\Y_{\nu,\e}}\nonumber\\
\le \|\oo J_{R,H}(\g_{m_\e})\|_{L_1(\O)}+ \oo C_{H,\nu,{m_\e}}\|\g\|_{\Y_{\nu,\e}},
\label{J3n}
\eaa
where $c_J$ is the same as in (\ref{J1}), and where \baaa
\oo C_{H,\nu,{m_\e}}\defi c_JC_H\sum_{k=2^{m_\e}+1}^\infty (2^{-k})^{\nu},
 \eaaa
 Clearly, $\oo C_{H,\nu,{m_\e}}$ is independent on $\g\in\Y_{\nu,\e}$, and, for any $c\in(1/2,1)$, $\oo C_{H,\nu,{m_\e}}$ is  bounded by a constant for all  $H\in (1/2,c), \e>0$.

Further, let \baaa
\xi^{(m_\e)}_k(t)=\rho(t,T^{(m_\e)}_k)=d_H\int_0^{T^{(m_\e)}_k}(t-s)^{H-3/2}dB(s).
\label{kkn}
\eaaa

Let $M_{\e}\defi \oo C_0^2 C_H^2\e_{m_\e}^{2H-1}$ and
\baaa
a_k\defi \int_{0}^{T}\|\g_{m_\e}(t)\|^2_{L_2(\O)}dt,\quad b_k\defi \int_{T^{(m_\e)}_k}^{T^{(m_\e)}_{k+1}}
\|\xi^{(m_\e)}_k(t)\|_{L_2(\O)}^2dt.
\eaaa
Clearly, 
\baaa \sum_{k=1}^n a_k=\int_{0}^{T}\|\g_{m_\e}(t)\|^2_{L_2(\O)}dt
\le \int_{0}^{T}\|\g(t)\|^2_{L_2(\O)}dt\le\|\g\|_{\Y_{\e,\nu}}^2. 
\label{ak}
\eaaa
As was shown for $\xi_k(t)$  in (\ref{kk3333}), we have that $b_k \le M_{\e}$ for all $k$. 

\par
We have that, for any $c\in (1/2,1)$,
\baa
&&\E|\oo J_{R,H}(\g_{m_\e})|\le\sum_{k=1}^{2^{m_\e}} \int_{T^{(m_\e)}_k}^{T^{(m_\e)}_{k+1}}\|\g_{m_\e}(t)\|_{L_2(\O)}\|\xi^{(m_\e)}_k(t)\|_{L_2(\O)}dt\le\sum_{k=1}^{2^{m_\e}} a_k^{1/2}b_k^{1/2}\nonumber \\&& \le \left( \sum_{k=1}^{2^{m_\e}} a_k\right)^{1/2}\left( \sum_{k=1}^{2^{m_\e}}  b_k\right)^{1/2}\le
 M_{\e}^{1/2} \cdot  2^{m_\e/2} \|\g\|_{\Y_{\e,\nu}}\le \w C\|\g\|_{\Y_{\e,\nu}}. \label{kk2n}
\eaa
for some $\w C=\w C(c,m_\e)>0$. We have used here the H\"older's inequality. 

It can be noted that the value
$m_\e$  in (\ref{kk2n}) is not increasing, since $\e>0$ is fixed.

By the definitions,  $\g_0(t)$ is $\GG_0$-measurable.
By the second estimate in Proposition \ref{propWR},
\baaa
\E|\oo J_{R,H}(\g_0)|\le\w C_0\|\g\|_{\L_{22}}\le \w C_0\|\g\|_{\Y_{\e,\nu}}.
\eaaa
for some $\w C_0=\w C_0(c)$.

The proof of Theorem \ref{ThM}(ii) follows from (\ref{J3n}) and (\ref{kk2n}). This completes the proof of Theorem \ref{ThM}.  $\Box$
 
{\em Proof of Proposition \ref{ThHolder}}
repeats the proof of  Theorem \ref{ThM}, given the adjustments
mentioned in the proof of Lemma  \ref{lemmaR}. $\Box$

The remaining part of the paper is devoted to the proof of Theorem \ref{ThL}.
We will use the notations from the proof of Theorem \ref{ThM} with the following amendment: since we consider variable $H\in[1/2,1)$, we include
corresponding $H$ as an index for a variable.

In particular, it follows from these notations that
\baaa
I_{W,H}(\g_n)=\sum_{k=1}^{n} P_{W,H,k}+I_{1/2}(\g_n),
\label{IHGw}\eaaa

It can be noted that
\baaa
d_H=\frac{H-1/2}{\G(H+1/2)}\to 0,\quad C_H=\frac{\sqrt{\G(H+1/2)^2 (H-1/2)}}{2-2H}\to 0\quad\hbox{as}\quad H\to 1/2+0.
\eaaa
\begin{lemma}\label{lemmaWw}  For any $\g_n\in\X_\e$, \baaa
\|I_{W,H}(\g_n)-I_{1/2}(\g_n)\|_{L_2(\O)}+\|\w I_{R,H}(\g_n)\|_{L_1(\O)}\to 0 \quad\hbox{as}\quad H\to 1/2+0
\label{II}
\eaaa
uniformly over any bounded in $\L_{22}$ set  of $\g_n\in \X_\e$.
\end{lemma}
\par
{\em Proof of Lemma \ref{lemmaWw}}.
For the operators $\G_k(\cdot,\cdot)=G_{H}(\cdot,T_{k+1},\cdot)$ introduced before  Lemma \ref{lemmaW}, we have that
 $\|\G_k(\cdot,g)\|_{L_2(T_k,T_{k+1})}\le\w c\|g\|_{L_2(T_k,T_{k+1})}$   for some $\w c>0$ that
is independent on $H\in(1/2,1)$. Similarly to the proof of Lemma \ref{lemmaW},
we have  that
\baaa
P_{W,H,k}
=\int_{T_{k-1}}^{T_k}  dB(\tau)[\G_{k-1}(\tau,\g_n)-\g_n(\tau)], \quad k=1,...,n.
\eaaa
These  integrals converge  in $L_2(\Omega,\GG_{T},\P)$.

Let 
\baaa
\a_{H,k}\defi  
\int_{T_{k-1}}^{T_k} |\G_{k-1}(\tau,\g_n)-\g_n(\tau)|^2d\tau.
\eaaa
We have that $\E \a_{H,k}=\E  P_{W,H,k}^2$ and 
\baaa
&&\E\|I_{W,H}(\g_n)-I_{W,1/2}(\g_n)\|_{L_2(\O)}^2=\E \left(\sum_{k=1}^n P_{W,H,k}\right)^2=\E\sum_{k=1}^n \a_{H,k}.
\eaaa

By   the properties of the Riemann--Liouville integral, we have  that
\baaa
\|\g_n-\G_k(\cdot,\g_n)\|_{L_2(T_{k-1},T_k)}\to 0
\quad\hbox{a.s. as}\quad H\to 1/2+0
\eaaa 
and 
\baaa
\|\G_k(\cdot,\g_n)\|_{L_2(T_{k-1},T_k)}\|\le\|\g_n\|_{L_2(T_{k-1},T_k)}\quad\hbox{a.s.}.
\eaaa
Hence  
\baaa
0\le \a_{H,k}\le \sqrt{2}\|\g_n\|_{L_2(T_{k-1},T_k)}\quad\hbox{a.s.}.
\eaaa
By Lebesgue's Dominated convergenceTheorem, it follows that
\baaa
\E\sum_{k=1}^n \a_{H,k}\to 0
\quad\hbox{a.s. as}\quad H\to 1/2+0.
\eaaa 
Hence
\baaa
\E\|I_{W,H}(\g_n)-I_{W,1/2}(\g_n)\|_{L_2(\O)}^2\to 0 \quad\hbox{as}\quad H\to 1/2+0.
\label{III}
\eaaa
\par
Further, we have that
\baaa
\E|\w I_{R,H}(\g_n)|\le \left(\E\int_{0}^{T}\g_n(t)^2dt\right)^{1/2} \left(\E\int_{0}^{T}\w\rho(t)^2dt\right)^{1/2}.
\eaaa
Similarly to the proof of Proposition \ref{propM}, we obtain that
\baa
\E\w\rho (t)^2&=&\int_{-\infty}^{\Ts}|f_t'(t,r)|^2dr=\frac{d_H^2}{2-2H}t^{2H-2}
\label{DR}
\eaa
and
\baaa
\E\int_{0}^T \w\rho (t)^2dt=\frac{d_H^2}{2(2-2H)}T^{2H-1}=\frac{c_hd_H}{4}T^{2H-1}\to 0 \quad\hbox{as}\quad H\to 1/2+0.
\label{EDR11} \eaaa
This completes the proof of Lemma \ref{lemmaWw}. $\Box$

\begin{lemma}\label{lemmaRw}  Let $\nu>0$, $\g\in \Y_{\nu,\e}$,
and $\{\g_n\}_{n=1}^\infty=\Z(\g)$.  In the notations introduced above, we have that
 \baaa
\|\oo J_{R,H}(\g_n)\|_{L_1(\O)}\to 0 \quad\hbox{as}\quad H\to 1/2+0
\label{IIn}
\eaaa
uniformly in $n>0$. \end{lemma}
\par
{\em Proof of Lemma \ref{lemmaRw}}.
Assume that $\g\in\Y_{\nu,\e}$ for some $\e>0$, and that $m_\e$ 
is defined by (\ref{m_e}). 
 It follows from equation (\ref{J2}) applied to $\oo J_{R}=\oo J_{R,H}$ that, for any and any $n>m_\e$,
\baaa
\E|\oo J_{R,H}(\g_n)|\le \|\oo J_{R,H}(\g_{m_\e})\|_{L_1(\O)}+ c_JC_H\sum_{k=m_\e+1}^n (2^{-k})^{\nu/2+H-1/2}\|\g\|_{\Y_{\nu,\e}}\nonumber\\
\le \|\oo J_{R,H}(\g_{m_\e})\|_{L_1(\O)}+ \oo C_{H,\nu,m_\e}\|\g\|_{\Y_{\nu,\e}},
\label{J33}
\eaaa
where $\oo C_{H,\nu,m_\e}$ is the same as in (\ref{J3n}); if $\nu>0$, then $\oo C_{H,\nu,m_\e}$ is  bounded by a constant for all  $H\in (1/2,1), \e>0$. In addition, we have that
\baaa
C_{H,\nu,m_\e}\to 0\quad \hbox{as}\quad H\to 1/2
\eaaa
uniformly in $n$. By  (\ref{kk2}),
$\|\oo J_{R,H}(\g_m)\|_{L_1(\O)}\to 0$ as  $H\to 1/2$.
This completes the proof of Lemma \ref{lemmaRw}. $\Box$

{\em Proof of Theorem \ref{ThL}}. Let $\g\in \Y_{\nu,\e}$ for any $\nu>0$  and $\e>0$.
Let  $\g_n=\Z(\g)$. We have to show that
$\E|I_H(\g)-I_{1/2}(\g)|\to 0$
as $H\to 1/2$. 
We have that
\baaa
\E|I_H(\g)-I_{1/2}(\g)|\le A_{1,H,n}+A_{2,H,n}+A_{3,n},
\eaaa
where
\baaa
A_{1,H,n}\defi\E|I_H(\g)-I_{H}(\g_n)|,\quad A_{2,H,n}\defi\E|I_H(\g_n)-I_{1/2}(\g_n)|,\quad A_{3,n}\defi\E|I_{1/2}(\g_n)-I_{1/2}(\g)|.
\eaaa
Clearly, $\|\g-\g_n\|_{\Y_{\nu,\e}}\to 0$ as $n\to +\infty$ for any $\e>0$. 

Let $c\in (1/2,1)$ be given. By Theorem \ref{ThM},
$A_{1,H,n}\to 0$  as $n\to +\infty$ uniformly in $H\in (1/2,c)$.
By Lemmata  \ref{lemmaWw}-\ref{lemmaRw},   $A_{2,H,n}\to 0$
as $H\to 1/2$ uniformly in $n$. Finally, by the properties of the It\^o integral, it follows that
$A_{3,n}\to 0$  as $n\to +\infty$. This completes  the proof of Theorem \ref{ThL}. $\Box$


\begin{thebibliography}{100}
\bibitem{Alos}
Al\'{o}s, E., Mazet, O., and Nualart, D. (2000).
 Stochastic calculus with respect to fractional Brownian
motion with Hurst parameter lesser than 1/2.
{\em Stochastic Processes and their Applications} 86 (2000) 121-139.
\bibitem{AN}
Al\'{o}s, E. and Nualart, D. (2003). Stochastic integration with respect to the fractional
Brownian motion. {\em Stoch. Stoch. Rep.} 75(3), 129-152.
\bibitem{BKiev}
Bender C., Sottinen T., Valkeila E. (2007).
Arbitrage with fractional Brownian motion? {\em Theory Stoch. Process.}, 13(1-2), 23-34 (Special Issue: Kiev Conference on Modern Stochastics).
\bibitem{B11}
Bender C., Sottinen T., Valkeila E.  (2011).
Fractional processes as models in stochastic finance. In: Di Nunno, Oksendal (Eds.), {\em  AMaMeF: Advanced Mathematical Methods for Finance}, Springer, 75-103.

\bibitem{B13}
Bender C. (2013).
An It\^o formula for generalized functionals of a
fractional Brownian motion with arbitrary Hurst
parameter. {\em Stochastic Processes and their Applications} 104, 81--106.
\bibitem{BPS}
Bender C., Pakkanen M.S., and Sayit H. (2015). Sticky continuous processes have consistent price
systems. {\em J. Appl. Probab.} 52, No. 2 , 586-594.
\bibitem{Ber}
Bertoin, J. (1989). Sur une int\'egrale pour les processus \'a $\a$
variation born\'ee, {\em The Annals of
Probability} 17, no. 4, 1521-1535.
\bibitem{Bi}
Biagini, F. and Oksendal, B. (2003).
Minimal variance hedging for fractional Brownian motion.  {\em Methods Appl. Anal.} 10 (3), 347-362.
\bibitem{Bj}
Bj\"ork T., Hult H. (2005). A note on Wick products and the fractional
Black-Scholes model.
{\em Finance and Stochastics} 9(2), 197-209.
\bibitem{Carm}
Carmona, P. Coutin, L., and Montseny, G. (2003). Stochastic integration with
respect to fractional Brownian motion. {\em Ann. Inst. H. Poincare Probab.
Statist.}, 39(1), 27-68.
\bibitem{CNS}
\c{C}etin U.,  Novikov A.,  Shiryaev A.N. (2013).
Bayesian sequential estimation of a drift of fractional Brownian motion.
 {\em Sequential Analysis: Design Methods and Applications }
32, Iss. 3, 288--296.\bibitem{C}
Cheridito, P. (2003). Arbitrage in fractional Brownian motion models. {\em Finance Stoch.} 7 (4),
533--553.
\bibitem{Ci}
  Ciesielski, Z., Kerkyacharian, G., and Roynette, B. (1992). Quelques espaces fonctionnels associ\'es \'a des processus gaussiens. {\em Studia Mathematica} 107 (1993), no. 2, 171-204.
\bibitem{DecU}  Decreusefond, L. and \"Ust\"unel, A.S.
(1999). Stochastic Analysis of the Fractional Brownian
Motion. {\em  Potential Analysis} 10 , no. 2, 177-214.
\bibitem{Dec0} Decreusefond, L. (2000). A Skohorod-Stratonovitch integral for the fractional Brownian motion,
Proceedings of the 7-th Workshop on stochastic analysis and related fields.
Birkhauser.
\bibitem{Dec}
Decreusefond, L. (2003).  Stochastic integration with respect to fractional Brownian motion.
In: P. Doukhan, G. Oppenheimer, M.S. Taqqu. (eds.), {\em Theory and Applications of Long-Range
Dependence}, Birkhauser, Boston, MA, pp. 203--226.

\bibitem{Duncan}
Duncan, T.E.,  Hu, Y., and Pasik-Duncan, B. (2000). Stochastic calculus for fractional Brownian
motion I, Theory. {\em SIAM Journal of Control and Optimisation} 38(2), 582?612.

\bibitem{Es}
Es-Sebaiya, K.,
    Ouassoub, I.,
 Oukninea, Y. (2009). Estimation of the drift of fractional Brownian motion. {\em Statistics \& Probability Letters}
79 (14), 1647--1653.
\bibitem{Fe}
 Feyel, D., and de La Pradelle, A. (1999). On Fractional Brownian Processes. {\em Potential Analysis}
10, no. 3, 273-288.
\bibitem{GN}
Gripenberg. G., and Norros, I. (1996). On the prediction of fractional Brownian motion.
{\em Journal of Applied Probability}
 33, No. 2, pp. 400-410.
\bibitem{G}
Guasoni, P. (2006). No arbitrage with transaction costs, with fractional Brownian
motion and beyond.  {\em Math. Finance} 16(2), 469--588.
\bibitem{H}
Hu Y. and Oksendal B. (2003). Fractional white noise calculus and applications to finance. {\em Infinite
dimensional analysis, quantum probability and related topics} 6(1), 1-32.
\bibitem{H2}
Hu Y., Oksendal B., Sulem A. (2003). Optimal consumption and portfolio in a Black-Scholes market driven by fractional Brownian motion.  {\em Infinite
dimensional analysis, quantum probability and related topics} 6(4), 519-536.
\bibitem{HZhou}
Hu Y., Zhou X.Y. (2005). Stochastic control for linear systems driven by fractional noises. {\em SIAM Journal on Control and Optimization} {\bf 43}, 2245-2277.
\bibitem{J}
Jumarie G. (2005). Merton's model of optimal portfolio in a Black-Scholes market driven by a fractional Brownian motion with short-range dependence,
{\em Insurance: Mathematics and Economics} 37, 585-598.
{\em Annals of the University of Bucharest (mathematical series)} 4 (LXII) (2013), 397?413

\bibitem{Ma}
Mandelbrot, B. B., Van Ness, J. W. (1968). Fractional Brownian motions, fractional noises and
applications. {\em SIAM Review} {\bf 10}, 422--437.
\bibitem{Mu}
Muravlev  A. A. (2013). Methods of sequential hypothesis testing for the drift of a fractional Brownian motion.
{\em Russ. Math. Surv.}  {\bf 68} (3), 577.
\bibitem{PT1}
Pipiras, V. and Taqqu, M.S. (2000). Integration questions related to fractional
Brownian motion. {\em Probab. Theor. Related Fields}, 118(2), 251-291.
\bibitem{PT2}
Pipiras, V. and Taqqu, M.S. (2001). Are classes of deterministic integrands
for fractional Brownian motion on an interval complete? {\em Bernoulli},
7(6), 873-897.
\bibitem{Pri}
Privault, N. (1998).  Skorohod stochastic integration with respect to non-adapted processes on
Wiener space. {\em Stochastics and Stochastics Reports} 65, 13--39.
\bibitem{R}
Rogers, L. C. G. (1997). Arbitrage with fractional brownian motion. {\em Mathematical Finance}
{\bf 7} (1),
95--105.
\bibitem{Shi}
Shiryaev A.N. (1998). On arbitrage and replication for fractal models. Research Report 30,
MaPhySto, Department of Mathematical Sciences, University of Aarhus.
\bibitem{Sal}
Salopek, D. M. (1998). Tolerance to arbitrage. {\em Stochastic Process. Appl.} 76 (2), 217-230.
\bibitem{Young}
Young, L.C. (1936).  An Inequality of H\"older Type, connected with Stieltjes integration. {\em  Acta
Math.} 67 , 251--282
\bibitem{Zahle}
 Z\"ahle, M. (1998). Integration with respect to fractal functions and stochastic calculus II.
 {\em Probability
Theory and Related Fields} 111 (3), pp. 333-374.
\end{thebibliography}
\end{document}